\documentclass[12pt]{amsart}
\usepackage{amssymb, amsfonts, amscd,amsmath}
\usepackage{calrsfs}

\newtheorem{theorem}{Theorem}[section]
\newtheorem{proposition}[theorem]{Proposition}
\newtheorem{corollary}[theorem]{Corollary}
\newtheorem{lemma}[theorem]{Lemma}

\newtheorem{remark}[theorem]{Remark}
\newtheorem{example}[theorem]{Example}

\newcommand{\cc}{\mathbb{C}}

\newcommand{\na}{\mathbb{N}}
\newcommand{\esv}{\mathfrak{M}V(X)}

\begin{document}


\title[Weighted algebras of holomorphic functions]{Spectra of weighted algebras of holomorphic functions}
\author{Daniel Carando \and Pablo Sevilla-Peris }

\address{Departamento de Matem\'atica, FCEN, Universidad de Buenos Aires\\
Intendente G\"uiraldes 2160, Ciudad Universitaria, C1428EGA\\
Buenos Aires, Argentina.}
\email{dcarando@dm.uba.ar}

\address{Departamento de Matem\'atica Aplicada, ETSMR, Universidad Polit\'ecnica de Valencia\\
Av. Blasco Ib\'a\~nez 21 46010\\
Valencia, Spain}
\email{Pablo.Sevilla@uv.es}
\subjclass[2000]{46G25, 46A45}

\keywords{weighted spaces and algebras, holomorphic functions,
spectrum, composition operators, algebra homomorphisms.
}

\thanks{The first author was partially supported by  PIP 5272 and PICT 05 17-33042. The second author was supported by the MECD Project MTM2005-08210 and grants GV-AEST06/092 and UPV-PAID-00-06}

\maketitle

\begin{abstract}
\noindent We consider weighted algebras of holomorphic functions
on a Banach space. We determine conditions on a family of weights
that assure that the corresponding weighted space is an algebra or
has polynomial Schauder decompositions. We study the spectra of
weighted algebras and endow them with an analytic structure. We
also deal with composition operators and algebra homomorphisms, in
particular to investigate how their induced mappings act on the
analytic structure of the spectrum. Moreover, a Banach-Stone type
question is addressed.
\end{abstract}

\section*{Introduction}

This work deals with weighted spaces of holomorphic functions on a
Banach space. If $X$ is a (finite or infinite dimensional) complex
Banach space and $U \subseteq X$ open and balanced, by a weight we
understand any continuous, bounded function $v:U \longrightarrow
[0,\infty[$. Weighted spaces of holomorphic functions defined by
countable families of weights were deeply studied by Bierstedt,
Bonet and Galbis in \cite{BiBoGa93} for open subsets of
$\mathbb{C}^{n}$ (see also \cite{BiBoTa98},\cite{Bo03}, \cite{BoDoLi99},\cite{BoDoLiTa98},\cite{BoFrJo03}). Garc{\'\i}a, Maestre and Rueda defined and
studied in
\cite{GaMaRu00} analogous spaces of functions defined on Banach spaces. We recall the definition of the weighted space%
\begin{multline*}
HV(U)= \{ f:U \rightarrow \mathbb{C} \text{ holomorphic }\colon \\
\|f\|_{v}=\sup_{x \in U} v(x) |f(x)| < \infty \text{ all } v \in V \}.
\end{multline*}
We endow $HV(U)$ with the Fr\'echet topology $\tau_{V}$ defined
by the seminorms $(\| \ \|_{v})_{v \in V}$. Since the family $V$
is countable, we can (and will throughout the article) assume it
to be increasing.

One of the most studied topics on weighted spaces of holomorphic
functions are the composition operators between them. These are
defined in a very natural way; if $\varphi : \tilde U \to U$ is a
holomorphic mapping and $V$, $W$ are two families of weights, the
associated composition operator $C_{\varphi} : HV(U) \to HW(\tilde
U)$ is defined as $C_{\varphi}(f)=f \circ \varphi$. There are a
number of papers on this topic, both in the finite dimensional and
infinite dimensional setting \cite{BoDoLi99}, \cite{BoDoLiTa98}, \cite{BoFr02},
\cite{BoFrJo03}, \cite{GaMaSe04}, \cite{GaMaSe06}. Among other things, the authors
study different properties of the operator $C_\varphi$ (when it is
well defined, continuous, compact, weakly compact or
completely continuous) in terms of properties of $\varphi$ (``size''
of its range, different kinds of continuity).

\smallskip
Our aim in this paper is to study the algebra structure of $HV(U)$
whenever it exists. We determine conditions on the family of weights
$V$ that are equivalent to $HV(U)$ being an algebra, and present
some examples. We also consider polynomial decompositions of
weighted spaces of holomorphic functions. For this, we give a
representation of the associated weight whenever the original weight
is radial. We show how the existence of a polynomial
$\infty$-Schauder decomposition and the presence of an algebra
structure are related, and how they lead us to the consideration of
weights with some exponential decay. Many of these results are new,
up to our knowledge, even for the several variables theory (i.e.,
for $X$ a finite dimensional Banach space). As an application of
these decompositions, we are able to present a somehow surprising
example: a reflexive infinite dimensional algebra of analytic
functions on $\ell_1$.

Whenever $HV(X)$ is an algebra, we study the structure of its
spectrum. For a symmetrically regular $X$ (see definitions in
Section~\ref{sec-espectro}), we endow the spectrum of $HV(X)$ with
a topology that makes it an analytic variety over $X^{**}$, much
in the spirit of Aron, Galindo, Garc\'{\i}a and Maestre's work for
the space of holomorphic functions of bounded type $H_b(X)$
\cite{ArGaGaMa96}. We show that any function $f\in HV(X)$ extends
naturally to an analytic function defined on the spectrum and this
extension can be seen to belong, in some sense, to $HV$.

We also study  algebra homomorphisms and composition operators
between spaces $HV(X)$ and $HV(Y)$, for $V$ a family of exponential
weights. Namely, we consider the algebra of holomorphic functions of
zero exponential type. This class of functions has been widely
studied in function theory in one or several variables since the
1930's \cite{Bo38,Bo54} and, even nowadays, its interest also arises
in areas such as harmonic and Fourier analysis, operator theory and
partial differential equations in complex domains. Every algebra
homomorphism induces a mapping between the respective spectra, and
we investigate how this induced mapping acts on the corresponding
analytic structures. We show that, contrary to the case for
holomorphic functions of bounded type \cite{CaGaMa05}, composition
operators induce mappings with good behaviour: they are continuous
for the analytic structure topology. We also characterize the
homomorphisms whose induced mappings are continuous. The results on
algebra homomorphisms allow us to address a Banach-Stone type
question. In this context, by a Banach-Stone question we mean the
following: if two Banach spaces have (algebraically and
topologically) isomorphic algebras of holomorphic functions, what
can be said about the spaces themselves? Some recent articles on
this kind of problems are \cite{CaGaMa05},\cite{Vi06}. A survey on
different types of Banach-Stone theorems can be found in
\cite{GarrJa02}. This question can be seen as a kind of converse of
the problem studied, for example, in \cite{DiDi98}, \cite{LaZa00}, \cite{CaCaGa00},
\cite{CaLa04}.

\smallskip
We now recall some definitions and fix some notation. We will
denote  duals by $X^{*}$ if $X$ is a Banach space and $E'$ if $E$
is a Fr\'echet space.

Given a weight $v$, its associated weight is defined as
\[
\tilde{v}(x)=\frac{1}{\sup\{|f(x)|:f\in Hv(U),\ \|f\|_{v} \leq
1\}} = \frac{1}{\| \delta_{x} \|_{(Hv(U))'}}\ ,
\]
where $\delta_{x}$ is the evaluation functional. It is a well known
fact \cite[Proposition 1.2]{BiBoTa98}, that $\|f\|_{v} \leq 1$ if
and only if $\|f\|_{\tilde{v}} \leq 1$ (hence $Hv(U)=H\tilde{v}(U)$
isometrically). We also have in \cite[Proposition 1.2]{BiBoTa98},
that $v \leq \tilde{v}$. However, it is not always true that there
exists a constant $C$ for which $\tilde{v} \leq C v$; the weights
satisfying this kind of equivalence with their associated weights
are called essential. A weight $v$ is called radial if
$v(x)=v(\lambda x)$ for every $\lambda \in \mathbb{C}$ with $|
\lambda |=1$ and norm-radial if $v(x_{1})=v(x_{2})$ whenever
$\|x_{1}\|=\|x_{2}\|$.

A set $A \subseteq U$ is called $U$-bounded if it is bounded and
$d(A,X \setminus U)>0$. Holomorphic functions of bounded type on $U$
are those that are bounded on $U$-bounded subsets. The space of all
these functions  is denoted by $H_{b}(U)$. By $H^{\infty}(U)$ we
denote the space of holomorphic functions that are bounded in $U$.
Following \cite[Definition 1]{GaMaRu00}, we say that a countable
family of weights $V$ satisfies Condition I if for every $U$-bounded
$A$ there is $v \in V$ such that $\inf_{x \in A} v(x)>0$. If $V$
satisfies Condition I, then $HV(U) \subseteq H_{b}(U)$ and the
topology $\tau_{V}$ is stronger than $\tau_{b}$ (topology of uniform
convergence on the $U$-bounded sets).

Given a Banach space $X$, the space of continuous, $n$-homogeneous
polynomials on $X$ is denoted by $\mathcal{P}(^{n} X)$. For a
given family of weights $V$, we write
$\mathcal{P}V(^{n}X)=\mathcal{P}(^{n} X) \cap HV(X)$.

A locally convex algebra will be an algebra $\mathcal{A}$ with a
locally convex structure given by a family of seminorms
$\mathcal{Q}$ so that for every $q \in \mathcal{Q}$ there exist
$q_{1}, q_{2} \in \mathcal{Q}$ and $C>0$ satisfying $q(xy) \leq C
q_{1}(x) q_{2}(y)$ (i.e., so that multiplication is continuous).
The spectrum of $\mathcal{A}$ is the space of non-zero continuous
multiplicative functionals. In the sequel, by ``algebra'' we will
mean a locally convex algebra.

We denote the spectrum of $H_{b}(X)$ by $\mathfrak{M}_{b}(X)$. Whenever $HV(X)$
is an algebra, we will denote its spectrum by $\mathfrak{M}V(X)$.

\section{Weighted algebras of holomorphic functions}

Next proposition determines conditions on the weights that make
$HV(X)$ an algebra. We thank our friend Jos\'e Bonet for helping
us fixing the proof, the final form of which is due to him.

\begin{proposition} \label{caract alg}
Let $U$ be an open and balanced subset of $X$ and $V$ be a family
of radial, bounded weights satisfying
Condition I. Then $HV(U)$ is an algebra if and only if for every
$v$ there exist $w \in V$ and $C>0$ so that
\begin{equation} \label{cond algebra}
v(x) \leq C \tilde{w}(x)^{2} \text{ for all } x \in U.
\end{equation}
\end{proposition}
\begin{proof}
Let us begin by assuming that $HV(U)$ is an algebra. Given $v \in V$
there are $C>0$ and $w_{1}, w_{2}$ so that $\|fg\|_{v} \leq C
\|f\|_{w_{1}} \|g\|_{w_{2}}$. Since $V$ is increasing, we can assume
$w_{1}=w_{2}=w$. Let us fix $x_{0} \in U$, and choose $f \in Hw(X)$
with $\|f\|_{w} \leq 1$ such that $f(x_{0}) = 1/ \tilde{w}(x_{0})$
(see \cite[Proposition 1.2]{BiBoTa98}). Taking the Ces\`aro means of
$f$ (see \cite[Section 1]{BiBoGa93}, or \cite[Proposition 4]{GaMaRu00})
we have a sequence $(h_{j})_{j}\subseteq HV(U)$ such that
$\|h_{j}\|_{w} \leq 1$ and $|h_{j}(x_{0})| \longrightarrow 1/
\tilde{w}(x_{0})$ as $j \to \infty$. We can assume that
$h_{j}(x_{0}) \neq 0$ for  $j$ large enough and we get
\begin{multline*}
v(x_{0}) = v(x_{0}) |h_{j}(x_{0})^{2}| \frac{1}{|h_{j}(x_{0})^{2}|}
\leq \| h_{j}^{2} \|_{v} \frac{1}{|h_{j}(x_{0})^{2}|} \\
\leq C \|h_{j} \|_{w}^{2} \frac{1}{|h_{j}(x_{0})|^{2}}
\leq C \frac{1}{|h_{j}(x_{0})|^{2}}.
\end{multline*}
Letting $j \to \infty$ we finally obtain \eqref{cond algebra}.
Conversely, if \eqref{cond algebra} holds, the fact that
$\|f\|_{w}= \|f\|_{\tilde{w}}$ for every $f$ easily gives that
$HV(U)$ is an algebra.
\end{proof}

The problem of establishing if a weighted space of functions is an algebra
was considered by L. Oubbi in \cite{Ou93} for weighted spaces of
continuous functions. In that setting, $CV(X)$ is an
algebra if and only if for every $v \in V$ there are $C>0$ and $w
\in V$ so that, for every $x \in X$
\begin{equation} \label{cond alg cont}
v(x) \leq C w(x)^{2}.
\end{equation}
Let us note that in our setting of holomorphic functions, since $w \leq \tilde{w}$,
if \eqref{cond alg cont} holds then $HV(U)$ is an algebra. On the other
hand, if the family $V$ consists of essential weights, then $HV(U)$ is
an algebra if and only if \eqref{cond alg cont} holds.

\medskip
Examples of families generating algebras can be constructed by
taking a weight $v$ and considering the family $V=\{ v^{1/n}
\}_{n=1}^{\infty}$. Since in the sequel we will want that these
families satisfy Condition I, we have to impose $v$ to be strictly
positive.

Not every weighted algebra can be constructed with ``$1/n$'' powers
of a strictly positive weight. In \cite[Example 14]{GaMaRu00},  a
family of weights $W=\{w_{n}\}_{n}$ satisfying Condition I so that
$H_{b}(U)=HW(U)$ is defined. If $U_{n}$ is a fundamental system of
$U$-bounded sets, each $w_{n}$ is defined to be $1$ on $U_{n}$ and
$0$ outside $U_{n+1}$ and such that $0\le w_n\le 1$. It is clear
that $w_{n} (x) \leq w_{n+1} (x)^{2}$ for every $x$. Let us see that
there is no positive weight $v$ such that $H_{b}(U)=HV(U)$ (where
$V$ is defined as before). We can view the identity $id:HW(U)
\longrightarrow HV(U)$ as a composition operator $id=C_{id_{U}}$;
then by \cite[Proposition 11]{GaMaSe06}  (see also \cite[Proposition 4.1]{BoFr02})
for each $n \in \mathbb{N}$ there exists $m$ so
that $C_{id_{U}} : H_{w_{m}}(U) \longrightarrow H_{v^{1/n}}(U)$ is
continuous. Then \cite[Proposition 2.3]{GaMaSe04}  (see also
\cite[Proposition 2.5]{BoDoLiTa98}) gives that $v(x)^{1/n} \leq K
\tilde{w}_{m}(x)$ for all $x$. Choose $x_{0} \not \in U_{m+1}$ and
we have $v(x_{0})=0$, so $v$ is not strictly positive.

Even if we drop the positivity condition on $v$ (or, equivalently,
Condition I on the family $V$), the fact that $v$ is zero outside
$U_{m+1}$ makes it easy to see that $H_b(U)$ cannot be $HV(U)$ if we
consider, for example, $U=X$ or $U=B_X$: take any holomorphic
function which is not of bounded type and dilate it so that it is
bounded on $U_{m+1}$.

Now we present some concrete examples of weighted algebras.

\begin{example} \label{H inf} \rm
Let $v$ be the weight on $B_X$ given by $v(x)=(1-\|x\|)^{\beta}$ and
let us define $V=\{ v^{1/n} \}_{n}$. Then,
$H^{\infty}(B_X)\varsubsetneq HV(B_X) \varsubsetneq H_b(B_X)$.

The first inclusion and the second strict inclusion are clear. To see
that the first one is also strict, we
choose $x^{*} \in X^{*}$ and $x_{0} \in X$ so that $\| x^{*} \|=
|x^{*} (x_{0})|=\| x_{0} \| = 1$ and $f(x)= \log(1- x^{*}(x))$.
Clearly $f$ is holomorphic and not bounded on the open unit ball
$B_{X}$. On the other hand, there exists a constant $C>0$ for
which
\[
(1-\|x\|)^{\beta} |\log(1-x^{*}(x))| \leq (1-\|x\|)^{\beta}
\log|1-x^{*}(x)|+C.
\]
Now, if $|1-x^{*}(x)|>1$, then $\log |1-x^{*}(x)| \leq 2$. If $|1-x^{*}(x)|<1$, then
$|1-x^{*}(x)| \geq \big|1- |x^{*}(x)| \big| \geq 1-\|x\|$ and
\[
(1-\|x\|)^{\beta} \log|1-x^{*}(x)| \leq (1-\|x\|)^{\beta}
\log(1-\|x\|).
\]
Since the mapping $t \in ]0,1] \rightsquigarrow (t^{\beta} \log t)$
goes to $0$ as $t$ does, we have $f \in HV(B_{X}) \setminus
H^{\infty}(B_{X})$.
\end{example}

\begin{example}\label{HV no es HB}\rm  Let $v$ be the weight on $X$ given by
$v(x)=e^{-\|x\|}$ and $V=\{ v^{1/n} \}_{n}$. When
$X=\mathbb{C}^{n}$, this weighted space $HV(\mathbb{C}^{n})$ is the
very well known algebra of entire functions of zero exponential type
(see, for example, \cite{Bo38}, \cite{Bo54}).

We have $H^{\infty}({X})\varsubsetneq HV(X) \varsubsetneq H_b(X)$.
To see that the second inclusion is strict, take $x^{*}\in X^{*}$
and define $f(x)=e^{x^{*}(x)^2}$. It is immediate that $f$ is a
holomorphic function of bounded type that is not in $HV(X)$.

On the other hand, $HV({X})$ cannot be $H^{\infty}({X})$.
\end{example}

We end this section by showing another example of a family that
gives an algebra but is not given by $\{v^{1/n}\}$. We thank our
friend Manolo Maestre for providing us with it.
\begin{example} \label{manolo} \rm
Let us consider a positive, decreasing function $\eta$ defined on $X$ and define
$v_{n}(x)=\sqrt[n]{\log \big(n(1+\|x\|) \big) \eta(\|x\|)}$. This clearly satisfies that
$v_{n}(x) \leq v_{2n}(x)^{2}$ for all $x$ but there is no $v$ such that $v_{n}=v^{1/n}$.
\end{example}

\section{Schauder decomposition and weighted algebras}

In this section, we consider two natural families of weights
obtained from a {decreasing continuous} function $\eta :
[0,\infty[ \longrightarrow ]0,\infty[$ such that $\lim_{t \to
\infty} t^{k} \eta(t)=0$ for every $k \in \mathbb{N}$. Let us
define two different families of weights,
$v_{n}(x)=\eta(\|x\|)^{1/n}$ and $w_{n}(x)=\eta
(\frac{\|x\|}{n})$, $n\in \mathbb N$.  Our aim is to study some
properties of the weighted spaces $HV(X)$ and $HW(X)$, where
$V=\{v_n\}_n$ and $W=\{w_n\}_n$. From what has already been said
in the previous section, $HV(X)$ is always an algebra. Note that
$v_1(x)=w_1(x)=\eta(\|x\|)$. For simplicity, we will write $v=v_1$
and $w=w_1$.\\

Following standard notation the real function $\eta$ can radially
extended to a weight on $\mathbb{C}$ by $\eta(z) = \eta(|z|)$ for $z
\in \mathbb{C}$ and its associated weight is given by
\[
\tilde{\eta}(t)=
\frac{1}{\sup \{|g(z)| \ : \ g \in H(\mathbb{C}) \ |g| \leq 1/\eta \mbox{ on } \mathbb{C} \}}.
\]
We then have weights on different spaces defined from the same
function $\eta$; it is natural now to ask how the associated
weights are related. The following proposition, showed to us by
Jos\'e Bonet, answers that question.

\begin{proposition} \label{eta asoc}
Let $X$ be a Banach space and $v$ a weight defined by
$v(x)=\eta(\|x\|)$ for $x \in X$. Then $\tilde{v}(x) = \tilde{\eta}(\|x\|)$
for all $x \in X$.
\end{proposition}
\begin{proof}
Let us fix $x \in X$ and choose $x^{*} \in X^{*}$ such that
$\|x^{*}\|=1$ and $x^{*}(x) = \|x\|$. If $h \in H(\mathbb{C})$ is
such that $|h| \leq 1/\eta$ then, for any $y \in X$,
\[
| (h \circ x^{*} )(y) | = |h(x^{*}(y))|  \leq
\frac{1}{\eta(x^{*}(y))} \leq \frac{1}{\eta(\|y\|)} =
\frac{1}{v(y)}.
\]
So we have  $\|h\circ x^*\|_v\le 1$ and hence
\begin{multline*}
\frac{1}{\tilde{\eta}(\|x\|)} = \sup \{|h(\|x\|)| \ : \ h \in H(\mathbb{C}), \ |h| \leq 1/\eta\} \\
= \sup \{|(h \circ x^{*} )(x)| \ : \ h \in H(\mathbb{C}), \ |h| \leq 1/\eta\} \\
\leq \sup \{|f(x)| \ : \ f \in H_v(X), \ \|h\|_v\le 1 \} =
\frac{1}{\tilde{v}(x)}
\end{multline*}
and $\tilde{v}(x) \leq \tilde{\eta}(\|x\|)$.

 Let us suppose now
that $\tilde{v}(x) < \tilde{\eta}(\|x\|)$ for some $x \neq 0$. Then
there exist $f \in H(X)$ with $\|f\|_v \leq 1$ such that $|f(x)| >
1/\tilde{\eta}(\|x\|)$. Let us define now $g :  \mathbb{C}
\rightarrow \mathbb{C}$ by $g(\lambda) = f(\lambda x /\|x\|)$;
clearly $g \in H(\mathbb{C})$ and $|g(\lambda)| \leq
1/\eta(\lambda)$ for all $\lambda \in \mathbb{C}$. Therefore
$|g(\|x\|)| \leq 1/\tilde{\eta}(\|x\|)$, but this contradicts the
fact that $g(\|x\|)=f(x)$. This gives that $\tilde{v}(x) = \tilde{\eta}(\|x\|)$ for every
$x \neq 0$. Both $\tilde{v}$ and $\tilde{\eta}$ are continuous since $\eta$ is so, then
we also have $\tilde{v}(0) = \tilde{\eta}(0)$
\end{proof}
\noindent As an immediate consequence of this result we have that
$v$ is essential if and only if $\eta$ is so.

\begin{remark} \label{tilde w}
\rm Proceeding as in the previous Proposition we can easily show
that $\tilde{w}_{n} (x) = \tilde{\eta} (\|x\|/n)$. Indeed, let us
consider $\mu(t)=\eta(t/n)$ for $t>0$. If $f \in H(X)$ is such that
$|f| \leq 1/\eta$, then the function defined by $g(x) = f(x/n)$ is
clearly holomorphic on $X$ and $|g| \leq 1/\mu$. From this,
$\tilde{w}_{n} (x)=\tilde{\mu} (\|x\|) \leq \tilde{\eta} (\|x\|/n)$.
On the other hand, suppose there is some $x_{0} \in X$ so that
$\tilde{\mu} (\|x_{0}\|) < \tilde{\eta} (\|x_{0}\|/n)$. We can find
$f \in H(X)$ such that $|f| \leq 1/\mu$ and $|f(x_{0})| >
1/\tilde{\eta(\|x_{0}\|/n)}$. Defining $h(x)=f(nx)$ we get the
desired contradiction.
\end{remark}

\medskip Our family $W$ was already defined and studied in
\cite[Example 16]{GaMaRu00}. By \cite[Theorem 11]{GaMaRu00},
$\big(\mathcal{P}W(^{n} X) \big)_{n}$ is an $\mathcal{S}$-absolute,
$\gamma$-complete decomposition of $HV(X)$ (see \cite[Definition 3.32]{LibroDi99}
and \cite[Definition 3.1]{Ka70}). Let us see that,
furthermore, it is an $\infty$-Schauder decomposition. Let us recall
that a Schauder decomposition $(F_{n})_{n}$ of a Fr\'echet space $F$
is an $R$-Schauder decomposition (\cite[Theorem 1]{GalMaRu00}),
whenever, for any $(x_{n})_{n}$ with $x_{n} \in F_n$,
$\sum_{n}x_{n}$ converges in $F$ if and only if $\limsup_{n}
\|x_{n}\|^{1/n}\le 1/R$. It is well known \cite[Lemma 6]{GalMaRu00}
that any $\infty$-Schauder decomposition is
$\mathcal{S}$-absolute.

By \cite[Example 16]{GaMaRu00}, $\mathcal{P}W(^{n} X) =
\mathcal{P}w(^{n} X)$ topologically for every $n$. However,
$\mathcal{P}w(^{n} X)$ is a Banach space with just one ``natural"
norm, namely $\|\cdot\|_{w}$, while $\mathcal{P}W(^{n} X)$ has many
possible norms. Since $\infty$-Schauder decompositions are sequences
of Banach spaces, we will always consider $\mathcal{P}W(^{n} X)$ as
a Banach space with the norm $\|\cdot\|_{w}$.

\begin{proposition} \label{infty Schauder}
$\big(\mathcal{P}W(^{n} X)\big)_{n}$ is an $\infty$-Schauder
decomposition of $HW(X)$.
\end{proposition}
\begin{proof}
We want to show that $\sum_{m} P_{m}$ converges in $\tau_{W}$ if and
only if $\lim_{m} \|P_{m}\|_{w}^{1/m}=0$.

Let us suppose first that
$\sum_{m} P_{m}$ converges in $\tau_{W}$. Taking a sequence
$\alpha_{m}=1$ for all $m$, since it is an $\mathcal{S}$-Schauder
decomposition,
\[
\|\sum_{m} P_{m} \|_{\alpha}  = \sum_{m} \|P_{m}\|_{w}
\]
converges. Then, given any $R>0$, we can take $n>R$ and
\begin{equation} \label{no se}
\begin{split}
\sup_{x \in X} &  |P_{m} (x)|\ \eta(\|x\|_{X}) R^{m}
\leq \sup_{x \in X} |P_{m} (x)| \  \eta(\|x\|_{X}) n^{m} \\
& = \sup_{x \in X} |P_{m} (n x) | \  \eta(\|y\|_{X})
= \sup_{y \in X} |P_{m} (y)| \  \eta(\frac{\|x\|_{X}}{n}) =
\|P_{m}\|_{w_{n}}.
\end{split}
\end{equation}
Hence $\sum_{m} \big(\sup_{x \in X} |P_{m} (x)| \
\eta(\|x\|_{X})\big) R^{m} < \infty$ for all $R>0$ and this implies
that $\lim_{m} \|P_{m}\|_{w}^{1/n}=0$.

Now, if $\lim_{m} \|P_{m}\|_{w}^{1/n}=0$, then  $\sum_{m}
\big(\sup_{x \in X} | P_{m} (x)| \ \eta(\|x\|_{X})\big) R^{m} <
\infty$ for all $R>0$. Using \eqref{no se}, $\sum_{m} \|P_{m}
\|_{w_n}$ converges for all $n$ and this completes the proof.
\end{proof}

The space $HW(X)$ is not necessarily an algebra. We want to find
now conditions on the weight that make $HW(X)$ an algebra and to
study how is  $HW(X)$ related to $HV(X)$ in this case.

\begin{proposition} \label{HW algebra}
$HW(X)$ is an algebra if and only if there exist $k>1$ and
$C>0$ so that,  for all $t$,
\begin{equation} \label{cond HW algebra}
\eta (kt) \leq C \tilde{\eta}(t)^{2}.
\end{equation}
If, furthermore, $\eta$ is essential, then $HW(X)$ is an algebra if
and only if there exist $k>1$ and $C>0$ so that,  for all $t$,
\begin{equation} \label{cond HW algebra ess}
\eta (kt) \leq C \eta(t)^{2}.
\end{equation}
In this case,  $HW(X) \hookrightarrow HV(X)$ continuously and there
exist positive constants $a,b$ and $ \alpha$ so that $\eta(t) \leq a
e^{-b t^{\alpha}}$ for all $t$.
\end{proposition}
\begin{proof}
By Proposition~\ref{caract alg} and Proposition \ref{eta asoc}, if $HW(X)$ is an algebra, given
$n=1$ there exist $C>0$ and $m$ such $\eta(t)\le C \tilde{\eta}(\frac
t{m})^{2}$ for all $t$. This clearly implies \eqref{cond HW algebra}. On
the other hand, if \eqref{cond HW algebra} holds, given $n$ we can
choose $m_{n}$ so that $m_{n}>kn$ and the fact that $\tilde{\eta}$ is
decreasing (because $\eta$ is decreasing \cite{BoDoLiTa98}), together with
Proposition~\ref{caract alg} and Remark \ref{tilde w}, give that
$HW(X)$ is an algebra.

Now, if $\eta$ is essential, condition \eqref{cond HW algebra} is
equivalent to \eqref{cond HW algebra ess}. In this case, $\eta(t)
\leq C^{2^n-1} \eta(t/k^{n})^{2^{n}}$ for all $t$ and $n$. Hence,
given $m$ let us take $n$ such that $2^{n} > m$, then since $\eta$ is decreasing,
\[
\left( \frac{\eta (t)}{\eta(0)} \right)^{1/m} \leq \left( \frac{\eta(t)}{\eta(0)} \right)^{1/2^{n}}
\leq C^{1-1/2^n} \frac{\eta(t/k^{n})}{\eta(0)^{1/2^{n}}} .
\]
This gives
\begin{equation}\label{cota eta}
\eta (t)^{1/m} \leq C^{1-1/2^n} \eta(0)^{1/m-1/2^{n}} \eta(t/k^{n}).
\end{equation}
This means that there is $K>0$ such that  $v_{m}(x) \leq K w_{k^n}(x)$ for all $x \in E$.
Therefore, $HW(X) \hookrightarrow HV(X)$ continuously.

Moreover, since
$\eta(t)\to 0$ as $t\to \infty$, we can choose $r$ such that
$C\eta(r)<1$. We have $\eta(k^{n}r) \leq C^{2^n-1}
\eta(r)^{2^{n}}\le (C\eta(r))^{2^n}$ for all $n$. Now, for any
$t>0$, let $n$ be such that $k^{n}r \leq t < k^{n+1}r$. We have
$$
\eta(t) \le \eta(k^{n}r) \le  (C\eta(r))^{2^n} \le
(C\eta(r))^{\frac12 (t/r)^{\log_{k}2}}$$ which is bounded by
$ae^{-bt^\alpha}$ for a proper choice of positive constants $a, b$
and $\alpha$.
\end{proof}

We have given conditions for $HW(X)$ to be an algebra. We also had
that $\big(\mathcal{P}W(^{n} X)\big)_{n}$ is an $\infty$-Schauder
decomposition of $HW(X)$. Knowing that the polynomials form a
Schauder decomposition of a space of holomorphic functions is
useful, since it allows to derive some properties of the space of
holomorphic functions (reflexivity, different approximation
properties, etc.) from the properties of the spaces of homogeneous
polynomials. The more we know about the decomposition (being it
absolute, complete, etc \dots), the more we can conclude about the
space itself. Let us check when the
polynomials are such a decomposition for $HV(X)$.

Let us first note that $\mathcal{P}V(^{n} X)=\mathcal{P}W(^{n}
X)=\mathcal{P}v(^{n} X)$. We consider in $\mathcal{P}V(^{n} X)$ the
norm $\|\cdot\|_v$. Then if $\big(\mathcal{P}V(^{n} X)\big)_{n}$ is
an $\infty$-Schauder decomposition of $HV(X)$, by \cite[Theorem 9]{GalMaRu00},
we get $HV(X)=HW(X)$. Since we know that $HW(X)$ always admits such a decomposition,
we have that the spaces of weighted
polynomials form an $\infty$-Schauder decomposition of $HV(X)$ if
and only if $HV(X)=HW(X)$. Moreover, we have

\begin{proposition} \label{HV infinito schauder}
If $\eta$ is essential, $\big(\mathcal{P}V(^{n} X)\big)_{n}$ is an $\infty$-Schauder
decomposition of $HV(X)$ if and only if $HV(X)=HW(X)$.

In this case, there exist positive constants
$a_1,a_1,b_1,b_2,\alpha_1$ and $ \alpha_2$ such that $$ a_1 e^{-b_1
t^{\alpha_1}} \le \eta(t) \le a_2 e^{-b_2 t^{\alpha_2}}$$ for all
$t$.
\end{proposition}
\begin{proof}
We only need to show the inequalities. If $HV(X)=HW(X)$, then $HW(X)$ is
an algebra and the second inequality follows from
Proposition~\ref{HW algebra}. On the other hand, if $HV(X)=HW(X)$ there must
exist $m\in\na$ and $C>0$ such that $w_2(x) \le C v_m(x) $ for all
$x\in X$. This means that
\[
\eta(t/2)\le C \eta(t)^{1/m}
\]
for all $t$. Now we can proceed as in the last part of the proof of
Proposition~\ref{HW algebra} to obtain the desired inequality.
\end{proof}

\begin{remark}\label{exponenciales} \rm
There is a whole class of functions $\eta$ for which $HV(X)$ and
$HW(X)$ coincide (and, then, they are algebras with a polynomial
$\infty$-Schauder decomposition). Indeed, for any $b, \alpha>0$ we
can define $\eta(t)=e^{-bt^\alpha}$. Since $\eta(t/n)\le
\eta(t)^{1/n^{[\alpha]}}$ and $\eta(t)^{1/n}\le
\eta(t/n^{1/[\alpha]})$, we have $HV(X)=HW(X)$ topologically.

On the other hand, Proposition~\ref{HV infinito schauder} shows that
any $\eta$ satisfying $HV(X)=HW(X)$ must be bounded below and above
by functions of this type.
\end{remark}

If we want $HV(X)$ to have a polynomial decomposition without being
$HW(X)$, we must then weaken our expectation on the type of
decomposition. The polynomials form an $\mathcal{S}$-Schauder,
$\gamma$-complete decomposition of the weighted space of holomorphic
functions  whenever the family is formed by norm radial weights
satisfying Conditions I and II' (see \cite[Theorem 11]{GaMaRu00}).
Condition I was already introduced. We say that a family of weights
satisfies Condition II' if for every $v$ in the family there exist
$C>0$, $R>1$ and $w$ in the family so that $v(x) \leq C w(Rx)$ for all $x$
\cite[Proposition 8]{GaMaRu00}. We can characterise Condition II'
in terms of the function $\eta$. Note that this condition also
imposes a relationship between $HV(X)$ and $HW(X)$

\begin{proposition} \label{HV CondII}
The family $V$ satisfies Condition II' if and only if there exist
$R>1$, and $\alpha, C >0$  so that,  for all $t$,
\begin{equation} \label{cond V CondII}
\eta (t)^{\alpha} \leq C \eta(Rt).
\end{equation}
In this case,  $HV(X) \hookrightarrow HW(X)$ continuously.
\end{proposition}
\begin{proof}
First of all, if $V$ satisfies Condition II', clearly given any $n$
there exist $m$, $R$ and $C$ so that $\eta (t)^{m/n} \leq C
\eta(Rt)$ for all $t$. On the other hand, if \eqref{cond V CondII}
holds, for any $n$ let us choose $m \geq \alpha n$. Then
\[
\frac{\eta(t)^{1/n}}{\eta(0)^{1/(\alpha n)}}
\leq C \left( \frac{\eta(Rt)}{\eta(0)} \right)^{1/(\alpha n)}
\leq C  \left( \frac{\eta (Rt)}{\eta(0)} \right)^{1/m}
\]
and this gives that Condition II' holds.

Now, if $V$ satisfies Condition II' then for any given $n$ and $k$
we have $\eta (t/n) \leq \eta (R^{k}t/n)^{\alpha^{k}}$. Let $k$ be
such that $R^{k}>n$ and $m$ such that $m-1 \leq 1/\alpha^{k} \leq
m$. The set $A=\{t : \eta(t) \geq 1\}$ is compact; let then $K=
\sup_{A} \eta(t)^{\frac{1}{1/\alpha^{k}}}/\eta(t)^{1/m}$ and we
have
\[
\eta(t/n) \leq \eta (\frac{R^{k}}{n}t)^{\alpha^{k}} \leq \eta(t)^{\frac{1}{1/\alpha^{k}}}
\leq K \eta(t)^{1/m}.
\]
This completes the proof.
\end{proof}

\bigskip
Suppose we have Banach spaces $Z$ and $X$ and a continuous dense
inclusion $Z \hookrightarrow X$ (in fact, any injective operator
would do, but for the sake of simplicity we will consider an
inclusion). If $\eta : [0,\infty[ \longrightarrow ]0,\infty[$ is
decreasing, we have the already studied families of weights given by
$v_{n}(x)=\eta(\|x\|_{X})^{1/n}$ and $w_{n}(x)=\eta
(\frac{\|x\|_{X}}{n})$. These can also be considered as weights on $Z$.
This allows to define the spaces $HV_X(Z)$ and $HW_X(Z)$. Since $Z$
is dense in $X$,  then $\sup_{x \in Z} w(x) |P(x)| = \sup_{x \in X}
w(x) |P(x)|$ for all $P$. Applying \cite[Theorem 9]{GalMaRu00} we
get that $HW_X(Z)=HW(X)$ topologically.

The following examples make use of this simple fact: if $v$ is the
weight on $X$ given by $v(x)=e^{-\|x\|}$, then $\mathcal
Pv(^nX)=\mathcal P(^nX)$ for all $n$ and
\begin{equation}\label{Pv=P} e^{-1}\|P\|\le \|P\|_v\le n^n
{e^{-n}}\|P\|\end{equation}
 (see also \cite[Example 16]{GaMaRu00}).

Although not the simplest one, our results allow us to give a straightforward example of a reflexive algebra of analytic
functions on $\ell_1$.
\begin{example} \label{Tsirelson}
\rm Let $T^*$ be the original Tsirelson
space and consider the natural inclusion $\ell_1\hookrightarrow
T^*$. Let $v(x)=e^{-\|x\|_{T^*}}$ (i.e., $\eta(t)=e^{-t}$). As in
Remark~\ref{exponenciales}, it is easy to see that
$HV_{T^*}(\ell_{1})=HW_{T^*}(\ell_{1})$ and $HV(T^*)=HW(T^{*})$. By
the comments above, we have that $HV_{T^*}(\ell_{1})$ and $HV(T^*)$
are isomorphic algebras. Moreover, for each $n$, $\mathcal{P}V(^{n}
T^*)=\mathcal{P}(^{n} T^*)$ and then $\mathcal{P}V(^{n} T^*)$ is
reflexive \cite{AlArDi84}. Since $(\mathcal{P}V(^{n} T^*))_n$ is an
$\infty$-Schauder decomposition of $HV(T^*)$, this algebra is
reflexive \cite[Theorem 8]{GalMaRu00}. Therefore,
$HV_{T^*}(\ell_{1})$ is a reflexive algebra. Note that any weight of
exponential type such as those presented in
Remark~\ref{exponenciales} would have worked, so we have a whole family of such reflexive algebras.
\end{example}

\begin{example} \label{ele2 1} \rm
Let $(a_{k})_{k} \in \ell_{2}$ be a sequence of positive scalars and
define $X_{a}=\{ x \colon \|x\|=\sum_{k} |x_{k}a_{k}| < \infty\}$.
We consider again $\eta(t)=e^{-t}$ and $v(t)=\eta(\|x\|_X)$. Hence,
$HV_{X_a}(\ell_{2})$ is an algebra and
$HV_{X_a}(\ell_{2})=HV(X_{a})$. Let us see what $HV_{X_a}(\ell_{2})$
looks like. More precisely, let us first identify the homogeneous
polynomials belonging to $HV_{X_a}(\ell_{2})$. Note that $X_{a}$ is
isometrically isomorphic to $\ell_{1}$, the isometry given by the
mapping $x \rightsquigarrow (a_{k} x_{k})_{k}$. Then, for each $n$
we have $\mathcal PV_{X_a}(^n\ell_{2})=\mathcal PV(^nX_{a})=\mathcal
P(^nX_{a})=\mathcal P(^n\ell_1)$, and the last isomorphism is the
composition operator associated to $x \rightsquigarrow
(a_{k}x_{k})_{k}$.

Write an $n$-homogeneous polynomial $P$ on $\ell_2$  as
\[
P(x)= \sum_{\overset{\alpha_{1},\dots , \alpha_{n}=1}{_{\alpha_{1} \geq \cdots \geq \alpha_{n}}}}^{\infty}
b_{\alpha} x_{\alpha_{1}} \cdots x_{\alpha_{n}}.
\]
Then $P \in \mathcal{P}V_{X_a}(^{n} \ell_{2})$ if and only if the
polynomial given by
\[
P_a(y)= \sum_{\overset{\alpha_{1},\dots , \alpha_{n}=1}{_{\alpha_{1}
\geq \cdots \geq \alpha_{n}}}}^{\infty} b_{\alpha}
\frac{y_{\alpha_{1}}}{a_{\alpha_1}} \cdots
\frac{y_{\alpha_{n}}}{{a_{\alpha_1}}}.
\]
belongs to $\mathcal{P}(^{n} \ell_{1})$. This happens if and only if
there exists $K>0$ such that
\begin{equation} \label{polin ele 2 1}
|b_{\alpha}| \leq K a_{\alpha_{1}} \cdots a_{\alpha_{n}}
\end{equation}
for all $\alpha$. This means that a polynomial belongs to $
\mathcal{P}V_{X_a}(^{n} \ell_{2})$ if and only if its coefficients
are controlled in some way by the sequence $(a_k)_k$. As a
particular case we have $\ell_{2}^{*} \cap HV(\ell_{2})=\{ (b_{k})
\colon |b_{k}| \leq K a_{k} \}$ (which coincides, of course, with
$X_{a}^{*}$). As in Example~\ref{HV no es HB}, we have  $HV(X_a)
\hookrightarrow H_{b}(X_{a})$ strictly.
\end{example}

\section{The spectrum} \label{sec-espectro}
Our aim is now to study the structure of the spectrum of $HV(X)$.
This is well known for the space of holomorphic functions of bounded
type, $H_{b}(X)$, when $X$ is symmetrically regular. A complex
Banach space $X$ is said to be (symmetrically) regular if every
continuous (symmetric) linear mapping $T:X \to X^{\ast}$ is weakly
compact. Recall that $T$ is symmetric if $Tx_1(x_2) = Tx_2(x_1)$ for
all $x_1,x_2 \in X$. The first steps towards the description of the
spectrum $\mathfrak{M}_{b}(X)$ of $H_b(X)$ were taken by Aron, Cole
and Gamelin in their influential article \cite{ArCoGa91}. In
\cite[Corollary 2.2]{ArGaGaMa96} Aron, Galindo, Garc\'{\i}a and
Maestre gave $\mathfrak{M}_{b}(U)$ a structure of Riemann analytic
manifold modeled on $X^{\ast\ast}$, for $U$ an open subset of $X$.
For the case $U=X$, $\mathfrak{M}_{b}(X)$ can be viewed as the
disjoint union of analytic copies of $X^{**}$, these copies being
the connected components of $\mathfrak{M}_{b}(X)$). In
\cite[Section 6.3]{LibroDi99}, there is an elegant exposition of all
these results. The study of the spectrum of the algebra of the space
of holomorphic functions of bounded type was continued in
\cite{CaGaMa05}. We continue in this trend by studying here
$\mathfrak{M}V(X)$. In this section we present the analytic
structure of $\mathfrak{M}V(X)$, in the spirit of the above
mentioned results.

If $f$ is a holomorphic function defined on a Banach space $X$, we denote by
$\bar{f}$ or $AB(f)$ the Aron-Berner extension of $f$ to $X^{**}$
(see \cite{ArBe78} and \cite{LibroDi99} for definitions and
properties).

The copies of $X^{**}$ are constructed in the following way: given
an element $\phi$ in the
spectrum, we lay a copy of $X^{**}$ around $\phi$ considering, for
each $z\in X^{**}$, the homomorphism that on $f\in HV(X)$ takes
the value $ \phi \big(x \in X \rightsquigarrow \bar{f}(x+z)\big)$.
If we let $z$ move in $X^{**}$, we obtain a subset of the spectrum
that is isomorphic to $X^{**}$. But this works only if $\phi$ can
act on the function $x \in X \rightsquigarrow \bar{f}(x+z)$, that
is, if this function belongs to $HV(X)$.

\begin{lemma} \label{suma}
Let $V$ be a family of weights satisfying Conditions I and II'
such that every $v$ is decreasing and norm radial; then the mapping $HV(X)
\longrightarrow HV(X)$ given by $f \rightsquigarrow f( \cdot + x)$
is well defined and continuous for every fixed $x \in X$.
\end{lemma}
\begin{proof}
The mapping in the statement can be viewed as a composition operator
$C_{\varphi_{x}}$, where $\varphi_{x}:X \longrightarrow X$ is given
by $\varphi_{x}(y)=x+y$. We use \cite[Proposition 11]{GaMaSe06}
(see also \cite[Proposition 4.1]{BoFr02}) to see that it is
continuous.

Since $V$ satisfies Condition II', given $v \in V$, we can take $R>1$ and
$w_{1}$ so that $v(y) \leq w_{1} (Ry)$ for all $y$. Then if
$\|y\| > \frac{1}{R-1}\|x\|$, then $\|x+y\| \leq R \|y\|$ and
\[
{v(y)}\le {w_{1}(Ry)} \le {w_{1}(x+y)} .
\]
Let now $w_{2}$ be so that $\inf_{\|y\| \leq \frac{1}{R-1}\|x\|}
w_{2}(y)=c_{1}>0$; then,
\[
\sup_{\|y\| \leq \frac{1}{R-1}\|x\|} \frac{v(y)}{w_{2}(y+x)}
<\infty.
\]
Choosing $w \geq \max(w_{1}, w_{2})$ we finally obtain for some
$K>0$,
\[
\sup_{y \in X} v(y) |f(x+y)| \leq \sup_{y \in X} \frac{v(y)}{w(x+y)} \sup_{y \in X} w(x+y) |f(x+y)| \leq K \|f\|_{w}.
\]
\end{proof}

Recall that we are considering $v$ to be a decreasing, norm radial
weight. Since $v$ is a function of the norm, we can consider it
defined both on $X$ and $X^{**}$.

Davie and Gamelin showed that the Aron-Berner extension is an
isometry for polynomials with the usual norm. They first prove a
more general version of this fact:  if $z \in X^{**}$, there is
$(x_{\alpha})_{\alpha} \subseteq X$ such that $\|x_\alpha\|\le
\|z\|$ for all $\alpha$ and $P(x_{\alpha}) \to \bar{P}(z)$ as
$\alpha \to \infty$, for all polynomial $P$ on $X$
\cite[Theorem 1]{DaGa89}. By using their result we show now that the Aron Berner
extension is also an isometry from $\mathcal P V(^nX)$ into $\mathcal
P V(^nX^{**})$.

If $P \in \mathcal{P}v(^{n}X)$, clearly $\|P\|_{v} \leq \| \bar{P}
\|_{v}$. Also we can choose $x_{\alpha}$ in such a way that
$\|x_{\alpha}\| \leq \|z\|$ and
\[
v(z) | \bar{P}(z)| \leq \lim_{\alpha} v(z) |P(x_{\alpha})| \leq \sup_{\alpha} v(x_{\alpha}) |P(x_{\alpha})| \leq \| P \|_{v}.
\]
Therefore,
\begin{equation}\label{AB es isometria}
\|P\|_{v} =  \| \bar{P} \|_{v}.
\end{equation}
This implies that the Aron-Berner
extension is a continuous homomorphism from $HV(X)$ in
$HV(X^{**})$. This was showed to us by M. Maestre in a more
general setting, namely if $v$ is continuous on straight lines or
$w^{*}$-continuous on spheres.

\medskip
In what follows we consider a positive decreasing function $\eta$
such that there is $C>0$ with
\begin{equation} \label{cond espectro}
\eta(s) \eta(t) \leq C \eta (s+t).
\end{equation}
A simple example of such a function is  $\eta(t)=e^{-t}$. We
consider the family of weights $v_{n}(x)=\eta(\|x\|_{X})^{1/n}$,
defined analogously on $X^{**}$. The space $HV(X)$ is an algebra
and, since \eqref{cond V CondII} in Proposition  \ref{HV CondII}
holds, $V$ satisfies Condition II' and the weighted polynomials form
a Schauder decomposition of $HV(X)$. Also, by \cite[Example 16]{GaMaRu00}
it contains all the homogeneous polynomials. In order to
study $\esv$ we follow the notation and trends of \cite[Section 6.3]{LibroDi99}
for $\mathfrak{M}_{b}(X)$. We reproduce the construction for the sake of completeness.

Linear functionals belong to $HV(X)$, so we can
define an onto mapping $\pi : \esv \longrightarrow X^{**}$ by $\pi
(\phi)=\phi |_{X^{*}}$. Since the Aron-Berner
extension is continuous, we can also define $\delta : X^{**}
\longrightarrow \esv$ given by $\delta(z) (f)= \bar{f}(z)$. For
any given $f \in HV(X)$ there is an associated mapping $f'' : \esv
\longrightarrow \mathbb{C}$ defined by $f''(\phi)=\phi(f)$. The
canonical embedding of $X$ into $X^{**}$ is denoted by $J_{X}$.

For a fixed $z \in X^{**}$, we consider $\tau_{z}(x)=J_{X}x+z$ for $x \in X$. Since there is no risk of confusion
we also denote $\tau_{z}:HV(X) \longrightarrow HV(X)$ the mapping given by
\[
(\tau_{z} f)(x) = \bar{f} (J_{X}x+z) = \bar{f}( \cdot +z) = (\bar{f} \circ \tau_{z})(x).
\]
By Lemma~\ref{suma} and the comments above on the Aron-Berner
extension this mapping is well defined. As a consequence, we get
$\phi \circ \tau_{z} \in \esv$ for every $\phi \in \esv$ and $z \in
X^{**}$. If $X$ is symmetrically regular, then
$\tau_{z+w}f=(\tau_{z} \circ \tau_{w})f$ for all $f\in H_b(X)$
\cite[Lemma 6.28]{LibroDi99}. Since $V$ satisfies Condition I, we
have $HV(X) \hookrightarrow H_{b}(X)$ and $\tau_{z+w}=\tau_{z} \circ
\tau_{w}$ on $HV(X)$.

Also, if $x^{*} \in X^{*}$, we have
$\tau_{z}(x^{*})=z(x^{*})+x^{*}$. For $\phi\in\esv$, $\phi(z(x^{*}))
=z(x^{*})$ and then,
\[
(\phi \circ \tau_{z} ) (x^{*})= \phi (z(x^{*})+x^{*}) = z(x^{*})+
\phi (x^{*}).
\]
In other words,  $\pi (\phi \circ \tau_{z})=\pi(\phi)+z$.

For any
pair $\phi \in \esv$ and $\varepsilon >0$ we consider
\[
V_{\phi, \varepsilon} = \{ \phi \circ \tau_{z} \colon z \in X^{**} \ , \ \|z\| < \varepsilon \}.
\]
As in \cite[Section 6.3]{LibroDi99} we obtain that
$\mathcal{V}_{\phi}=\{V_{\phi, \varepsilon}\}_{\varepsilon >0}$ is a
neighbourhood basis at $\phi$ for a Hausdorff topology on $\esv$
whenever $X$ is symmetrically regular. Moreover,
$\pi(\phi)=\pi(\psi)$ if and only if $\phi=\psi$ or $V_{\phi, r}
\cap V_{\psi,s}= \O$ for all $r,s$; also $\esv$ is a Riemann domain
over $X^{**}$ whose connected components are ``copies'' of $X^{**}$.

As we have already mentioned, Condition I assures that $HV(X)
\hookrightarrow H_{b}(X)$. Moreover, all the polynomials belong to
$HV(X)$, so the inclusion has dense range. Hence, we have a one to one
identification $\mathfrak{M}_{b}(X) \hookrightarrow \esv$. We do not
know whether or not they are equal. Note that they both consist of
``copies'' of $X^{**}$.

We have the following commutative diagram\\
\begin{center}
\begin{picture}(130,90)
\put(0,15){$X$}
\put(0,85){$X^{**}$}
\put(110,85){$\esv$}
\put(120,15){$\mathbb{C}$}
\put(12,88){\vector(1,0){90}}
\put(12,18){\vector(1,0){105}}
\put(123,83){\vector(0,-1){57}}
\put(6,25){\vector(0,1){57}}
\put(-8,50){$J_{X}$}
\put(62,90){$\delta$}
\put(126,50){$f''$}
\put(62,5){$f$}
\put(10,83){\vector(2,-1){110}}
\put(65,60){$\bar{f}$}
\end{picture}
\end{center}

In the case of $H_{b}(X)$, the function $f''$ is holomorphic on
$\mathfrak{M}_{b}(X)$ and is, in some sense,  of bounded type. We
show now that something analogous happens in our situation. By the
Riemann domain structure of $\esv$, ``holomorphic'' means that $f''
\circ (\pi |_{V_{\phi, \infty}})^{-1}$ is holomorphic on $X^{**}$
for all $\phi \in \esv$, where $V_{\phi, \infty} =
\bigcup_{\varepsilon >0} V_{\phi, \varepsilon}$.

Given a weight $v$ defined on $X$, we define the corresponding weighted norm for $n$-linear mappings:
\[
\|A\|_{v} = \sup_{x_{1}, \dots , x_{n} \in X} |A(x_{1}, \dots , x_{n})| \  v(x_{1}) \cdots v(x_{n}).
\]
If $P \in \mathcal{P}(^{n} X)$, we denote the associated symmetric
$n$-linear mapping by $\check{P}$. For   a symmetric $n$-linear
mapping $A$, by $A(x^{k}, y^{n-k})$ we mean the mapping $A$ acting
$k$-times on $x$ and $(n-k)$ times on $y$.

\begin{lemma}  \label{multilineal}
Let $\eta$ be a positive, decreasing function satisfying \eqref{cond espectro} and $v(x)=\eta(\|x\|)$. Then,
for any $P \in \mathcal{P}v(^{n}X)$,
\[
\|\check{P}\|_{v} \leq \frac{C^{n}}{n!} \|P\|_{v}
\]
where $C$ is the constant in \eqref{cond espectro}.
\end{lemma}
\begin{proof}
For any choice of $x_{1}, \dots , x_{n} \in X$ we have,  using
\eqref{cond espectro} and the polarization formula,
\[
\begin{split}
|\check{P} (x_{1}, \dots , x_{n})| &  \  v(x_{1}) \cdots v(x_{n}) \\
\leq \frac{1}{2^{n}n!} \sum_{\varepsilon_{i} = \pm 1} &  P(\varepsilon_{1} x_{1} + \cdots + \varepsilon_{n} x_{n})
\  v(x_{1}) \cdots v(x_{n})
\leq \frac{C^{n}}{n!} \|P\|_{v}.
\end{split}
\]
\end{proof}

The following result is analogous to \cite[Proposition 6.30]{LibroDi99} and follows the same steps.
\begin{theorem} \label{espectro}
Let $X$ be symmetrically regular and $\eta$ be a positive, decreasing function satisfying \eqref{cond espectro}. Let
$V$ be defined by $v_{n}(x)=\eta(\|x\|)^{1/n}$. Then, for every $f \in HV(X)$, the associated function
$f'' : \esv \longrightarrow \mathbb{C}$ given by $f''(\phi) = \phi(f)$ is holomorphic.
\end{theorem}
\begin{proof}
For any  $\phi \in \esv$ and $z \in X^{**}$ we have
\[
\big(f'' \circ (\pi |_{V_{\phi, \infty}})^{-1}\big) (\pi(\phi)+z) = f''(\phi \circ \tau_{z})
= (\phi \circ \tau_{z})(f) = \phi(\tau_{z}f).
\]
Hence we need to prove that the mapping $z \in X^{**}
\rightsquigarrow \phi(\tau_{z}f)= \phi \big( x \mapsto \bar{f} (J_{X}x+z)\big)$ is holomorphic.

 Let us consider the polynomial expansion at zero: $f=\sum_{n}
P_{n}$, where $P_{n} \in \mathcal{P}(^{n} X)$ for all $n$. What we
need then is to show that the function $z \rightsquigarrow \phi
\big( x \mapsto \sum_{n} \bar{P}_{n}(z)(x)\big)$ is holomorphic.
To see it, this sum must converge for the topology $\tau_{V}$.
We write $A_{n}=\check{P }_{n}$. For $z \in X^{**}$ and $0 \leq k
\leq n$ define $P_{n,k,z}:X \longrightarrow \mathbb{C}$ by
$P_{n,k,z}(x)=\bar{A}_{n}(J_{X} x^{n-k}, z^{k})$; this is clearly
an $(n-k)$-homogeneous polynomial. Let us see that $P_{n,k,z}$
belongs to $\mathcal{P}V(^{n-k}X)$. For any $v \in V$, we set
$w_1=v^{1/(n-k)}$ and $w_2=v^{1/k}$. Then, choosing $w \geq \max
(w_{1},w_{2})$ we get
\begin{eqnarray*}
\|P_{n,k,z}\|_{v} &=& \sup_{x \in X} |\bar{A}_{n}(J_{X} x^{n-k}, z^{k})| v(x) \\
&\leq& \sup_{x \in X} |\bar{A}_{n}(J_{X} x^{n-k}, z^{k})|
\big(v(x)^{1/(n-k)} \big)^{n-k} \frac{1}{v(z)}
\big(v(z)^{1/k} \big)^{k} \\
&=& \sup_{x \in X} |\bar{A}_{n}(J_{X} x^{n-k}, z^{k})|
w_{1}(x)^{n-k} w_{2}(z)^{k} \frac{1}{v(z)} \\ &\le&
\|\bar{A}_{n}\|_{w} \frac{1}{v(z)}.
\end{eqnarray*}
Now we apply Lemma \ref{multilineal} to obtain
\begin{equation}\label{pnkz}
\|P_{n,k,z}\|_{v} \leq \|\bar{A}_{n}\|_{w} \frac{1}{v(z)} \leq
\frac{1}{v(z)} \frac{C^{n}}{n!} \|\bar{P}_{n}\|_{w} = \frac{1}{v(z)}
\frac{C^{n}}{n!} \|P_{n}\|_{w}.
\end{equation}

Proceeding as in \cite[Section 6.3]{LibroDi99}:
\begin{align*}
(\tau_{z}f)(x)& =\bar{f}(J_{X}x+z)= \sum_{n=0}^{\infty} \bar{P}(J_{X}x+z) \\
& = \sum_{n=0}^{\infty} \left( \sum_{k=0}^{n}
\binom n k \bar{A}_{n}(J_{X}x^{n-k},z^{k}) \right) \\
& = \sum_{n=0}^{\infty} \left( \sum_{k=0}^{n} \binom n k P_{n,k,z}
\right) (x).
\end{align*}
This gives a pointwise representation of the function. This series converges in
$\tau_{V}$; indeed if $v \in V$, inequality \eqref{pnkz} gives
\begin{align*}
\sum_{n=0}^{\infty} \sum_{k=0}^{n} \binom n k \sup_{x \in X} v(x)
|P_{n,k,z} (x)| & \leq \sum_{n=0}^{\infty} \sum_{k=0}^{n} \binom n k
\frac{C^{n}}{n!} \frac{1}{v(z)} \|P_{n}\|_{w} \\
& \leq K \frac{1}{v(z)} \sum_{n=0}^{\infty} \|P_{n}\|_{w}.
\end{align*}
Since $\eta$ is strictly positive, so is $v$ and by \cite[Lemma 10]{GaMaRu00}
the last series converges. Hence, for each $z \in X^{**}$,
the series $\sum_{n=0}^{\infty} \sum_{k=0}^{n} \binom n k
 P_{n,k,z}$ converges in $\tau_{V}$ to $\tau_{z}f$.
Then we can write
\[
\phi (\tau_{z}f) = \sum_{n=0}^{\infty} \sum_{k=0}^{n} \binom n k
\phi( P_{n,k,z}).
\]
Let us consider now the $k$-homogeneous polynomial $P_{n,k}:z \in X^{**}
\longrightarrow \phi(P_{n,k,z})$ and see that it is continuous. We
fix $w_{\phi} \in V$ such that $|\phi(h)| \leq M \|h\|_{w_{\phi}}$
for all $h \in HV(X)$. Note that $w_{\phi}$ coincides with
$\eta(\|\cdot\|)^{1/r}$ for some $r$. Let $z\in B_{X^{**}}$, by
\eqref{pnkz},
\[
|\phi(P_{n,k,z})| \leq M \|P_{n,k,z}\|_{w_{\phi}} \leq M
\frac{C^{n}}{n!} \|P_{n}\|_{w_{\phi}} \frac{1}{w_{\phi}(z)} \leq M
\frac{1}{n!} \|P_{n}\|_{w_{\phi}} \frac{1}{\eta(1)^{1/r}}.
\]
This means that $P_{n,k}$ is bounded and therefore
$Q_{n}=\sum_{k=0}^{n} \binom n k \phi( P_{n,k,z}) \in
\mathcal{P}(^{n} X^{**})$. Since
$\phi(\tau_{z}f)=\sum_{n=0}^{\infty} Q_{n}(z)$,  $\phi(\tau_{z}f)$
is a holomorphic function of $z$.
\end{proof}

We have shown that $f'' \in H(\esv)$. We can even get that in some
sense it ``belongs to $HV(\esv)$''. Let $\phi\in\esv$ and choose
$w_\phi$ as before. For any $v \in V$, let $u \geq
\max(w_{\phi},v)$. We have
\begin{multline*}
|f''(\phi\circ\tau_z)|v(z) \leq
\sum_{n=0}^{\infty}
\sum_{k=0}^{n} \binom n k \lvert \phi (P_{n,k,z})\rvert \ v(z) \\
\leq \sum_{n=0}^{\infty} \sum_{k=0}^{n} \binom n k M \| P_{n,k,z} \|_{w_{\phi}} v(z)
\leq  \sum_{n=0}^{\infty} \sum_{k=0}^{n} \binom n k M \| P_{n,k,z} \|_{u} u(z) \\
\leq \sum_{n=0}^{\infty} \sum_{k=0}^{n} \binom n k
M \frac{C^{n}}{n!} \frac{1}{u(z)} \| P_{n} \|_{u} u(z)
\leq M K \sum_{n=0}^{\infty} \| P_{n} \|_{u}.
\end{multline*}
which is a finite constant by \cite[Lemma 10]{GaMaRu00}. Therefore,
$f''$ belongs to $HV$ of each copy of $X^{**}$ in the spectrum.

\section{Algebra homomorphisms between weighted algebras}
We now consider the weight $v(\cdot)=e^{-\|\cdot\|}$ defined on any
Banach space, and the associated family $V=\{v^{1/n}\}_{n}$. This
family is given by $\eta(t)=e^{-t}$ and obviously satisfies
\eqref{cond espectro}. Moreover, $V$ and $W$ coincide, and
consequently the weighted spaces of polynomials are an
$\infty$-Schauder decomposition of the algebra $HV(X)$ for any
Banach space $X$. As mentioned above, $HV(X)$ is the algebra of
holomorphic functions of zero exponential type on $X$.

We now study continuous algebra homomorphisms
$A:HV(X) \longrightarrow HV(Y)$ and start by considering
composition operators.

First, just a remark: if $f$ is a holomorphic function such that
there exist $A,B>0$ with $|f(y)| \leq A \|y\|+B$ for all $y \in Y$,
then by the Cauchy inequalities,
\[
\left\| \frac{d^{k}f(0)}{k!} \right\|_{B_{Y}} \leq \frac{1}{r^{k}}
\|f\|_{rB_{Y}} \leq \frac{Ar+B}{r^{k}}.
\]
Unless $k=0$ or $k=1$, this goes to $0$ as $r$ goes to $\infty$.
Hence $f$ is affine: there exist $y^{*} \in Y^{*}$ and $C>0$ so that
$f(y)=y^{*}(y)+C$.

\begin{lemma} \label{serafin}
Let $A:HV(X) \longrightarrow HV(Y)$ be an algebra homomorphism. Then
$Ax^{*}$ is a degree 1 polynomial for all $x^{*} \in X^{*}$ (i.e.
$A$ maps linear forms on $X$ to affine forms on $Y$).
\end{lemma}
\begin{proof}
Since $A$ is continuous, given $n$, there exist $m$ and $C>0$ so that, for every $f \in  HV(X)$
\[
\sup_{y \in Y} |Af(y)| \ e^{-\frac{\|y\|}{n}} \leq C \sup_{x \in X} |f(x)| \  e^{-\frac{\|x\|}{m}}.
\]
Let us take $x^{*} \in X^{*}$ and define $f(x)=\sum_{x=0}^{M} \frac{x^{*}(x)^{j}}{\|x^{*}\|^{j} m^{j} j!}
\in HV(X)$. Since $A$ is an algebra homomorphism
\begin{multline*}
\sup_{y \in Y} \left| \sum_{j=0}^{M} \frac{(Ax^{*})(y)^{j}}{\|x^{*}\|^{j} m^{j} j!} \right| \ e^{-\frac{\|y\|}{n}}
\leq C \sup_{x \in X} \left| \sum_{j=0}^{M} \frac{x^{*}(x)^{j}}{\|x^{*}\|^{j} m^{j} j!} \right| \ e^{-\frac{\|x\|}{m}} \\
\leq C \sup_{x \in X} \sum_{j=0}^{M} \frac{|x^{j}|}{m^{j} j!} \ e^{-\frac{\|x\|}{m}}
\leq C \sup_{x \in X} e^{\frac{\|x\|}{m}} \ e^{-\frac{\|x\|}{m}} =C.
\end{multline*}
This holds for every $M$; hence $\sup_{y \in Y} \left|
e^{\frac{Ax^{*}(y)}{\|x^{*}\| m}} \right|\ e^{-\frac{\|y\|}{n}} \leq
C$. Then $\Re (\frac{Ax^{*}}{\|x^{*}\|} (y)) \leq K_{1} \|y\| + K_{2}$ for
all $y \in Y$. Also, if $| \lambda | =1$ we have $\Re (\lambda
\frac{Ax^{*}}{\|x^{*}\|} (y)) = \Re ( A \frac{\lambda x^{*}}{\|x^{*}\|} (y)) \leq
K_{1} \|y\| + K_{2}$. This gives $\left| A
\frac{x^{*}}{\|x^{*}\|}(y) \right| \leq K_{1} \|y\|+K_{2}$ for all
$y \in Y$. But this implies that $A \frac{x^{*}}{\|x^{*}\|}$ is
affine on $y$; hence so is $Ax^{*}$.
\end{proof}

\begin{corollary} \label{composicion}
If the composition operator $C_{\varphi}:HV(X) \longrightarrow HV(Y)$ is continuous, then $\varphi$ is affine.
\end{corollary}
\begin{proof}
By Lemma \ref{serafin}, $x^*\circ\varphi=C_{\varphi}(x^{*})$ is
affine. Since weakly affine mappings are affine, we obtain the
conclusion.
\end{proof}

It is clear that Lemma~\ref{serafin} and Corollary~\ref{composicion}
are not valid for operators from $H_b(X)$ to $H_b(Y)$. Indeed, for
any $\varphi\in H_b(Y,X)$, the composition operator $C_{\varphi}$ is
well defined and continuous from $H_b(X)$ to $H_b(Y)$. In some
cases, one may even obtain a non-affine bianalytic $\varphi$.
Indeed, if $f$ is any entire function on $\cc$, the Henon mapping
$h:\cc^2\to \cc^2$ given by $h(z,u):=(f(z)-cu,z)$ is bianalytic and,
of course, is not affine unless $f$ is. Henon-type mappings in
infinite dimensional Banach spaces were used in \cite[Theorem 35]{CaGaMa05}
to obtain homomorphisms with particular behaviour. See
comments below, after Corollary~\ref{abcomposicion}.

\medskip

As an application of the previous results, we obtain a
Banach-Stone type theorem for $HV$.

\begin{theorem} \label{Banach Stone}
If $HV(X) \cong HV(Y)$ as topological algebras, then $X^{*} \cong Y^{*}$.\\
If moreover both $X$ and $Y$ are symmetrically regular or $X$ is regular, then $HV(X)
\cong HV(Y)$ if and only if $X^{*} \cong Y^{*}$.
\end{theorem}
\begin{proof}
Let $A:HV(X) \longrightarrow HV(Y)$ be an isomorphism; by
Lemma~\ref{serafin}, $Ax^{*}$ is affine for every $x^{*} \in
X^{*}$. Let us define $S:X^{*} \longrightarrow Y^{*}$ by
$Sx^{*}=Ax^{*}-Ax^{*}(0_{Y})$. Clearly, $S$ is linear and
continuous. We consider also $\tilde{S}:Y^{*} \longrightarrow
X^{*}$ given by $\tilde{S}y^{*}= A^{-1} y^{*}
-A^{-1}y^{*}(0_{X})$. Taking into account that $Ax^{*}(0_{Y})$ and
$A^{-1}y^{*}(0_{X})$ are constants and that constants are
invariant for both $A$ and $A^{-1}$, it is easily seen than $S$
and $\tilde{S}$ are inverse one to each other. So $X^*$ and $Y^*$
are isomorphic.

If $X$ and $Y$ are symmetrically regular and $S:X^{*}
\longrightarrow Y^{*}$ is an isomorphism, by \cite[Theorem 4]{LaZa00}
the mapping $\hat{S}:\mathcal{P}(^{n} X) \longrightarrow
\mathcal{P}(^{n} Y)$ given by $\hat S(P)=\bar P\circ S^*\circ J_Y$
is an isomorphism. Since   $\mathcal{P}(^{n} X )$ and
$\mathcal{P}(^{n} Y)$ coincide with $\mathcal{P}v(^{n} X)$ and
$\mathcal{P}v(^{n} Y)$, we have that $\hat S$ is an isomorphism
between the weighted spaces of polynomials. However, we need an
estimation of the norm of $\hat S$ as an operator between
$\mathcal{P}v(^{n} X)$ and $ \mathcal{P}v(^{n} Y)$ to obtain the
isomorphism between the algebras. Since $v$ is decreasing and by
\eqref{AB es isometria} the Aron-Berner extension is an isometry
between the weighted spaces of polynomials we have
\begin{eqnarray*}
\sup_{y\in Y} v(y)|\hat S(P)(y)|& =& \sup_{y\in Y} v(y)|\bar P(
S^*(J_Y(y)))| \\ &=& \|S\|^n \sup_{y\in Y} v(y)\left|\bar P\left(\frac{ S^*(J_Y(y))}{\|S\|}\right)\right|\\
&=& \|S\|^n \sup_{y\in Y} v\left(\frac{ S^*(J_Y(y))}{\|S\|}\right)
\left|\bar P\left(\frac{ S^*(J_Y(y))}{\|S\|}\right)\right|\\ &\le&
\|S\|^n\|P\|_v
\end{eqnarray*}
Hence $\|\hat S(P)\|_v\le \|S\|^n\|P\|_v$ and analogously for
$\hat{S}^{-1}$. The fact that $\mathcal{P}v(^{n} X)$ and
$\mathcal{P}v(^{n} Y)$ are respectively $\infty$-Schauder
decompositions of $HV(X)$ and $HV(Y)$, \cite[Theorem 1]{GalMaRu00}
and the multiplicative nature  of the Aron-Berner extension give the
conclusion.

If either $X$ or $Y$ are regular, we proceed analogously using
\cite[Theorem 1]{CaCaGa00}.
\end{proof}

\smallskip
The spectrum of $HV(X)$ is formed by a number of copies of
$X^{**}$ and each one of them is a connected component of $\esv$.
This can be viewed as if each copy of $X^{**}$ were a ``sheet'' and
all those ``sheets'' were laying one over the other in such a way
that all the points in a vertical line are projected by $\pi$ on
the same element of $X^{**}$.

Every algebra homomorphism $A:HV(X) \longrightarrow HV(Y)$ induces
a mapping $\theta_{A}:\mathfrak{M}V(Y) \longrightarrow \esv$
defined by $\theta_{A}(\phi)=\phi \circ A$. The sheets (copies of
$Y^{**}$) are the connected components of $\mathfrak{M}V(Y)$. By
the analytic structure of $\mathfrak{M}V(Y)$, $\theta_{A}$ is
continuous if and only if $\theta_{A}$ maps sheets into sheets. We
want to characterize the continuity of $\theta_A$. In order to
keep things simple and readable we change slightly our notation.
From now on the elements of the biduals will be denoted by
$x^{**}$ and $y^{**}$. Also, we will identify $X^{**}$ and
$Y^{**}$ with their images $\delta(X)$ and $\delta(Y)$ in the
respective spectra.

\begin{theorem} \label{final}
Let $X$ and $Y$ be symmetrically regular Banach spaces and $A:HV(X)
\longrightarrow HV(Y)$ an algebra homomorphism. Then, the following
are equivalent.\\
$(i)$ \  There exist $\phi \in \esv$ and $T:Y^{**} \longrightarrow
X^{**}$ affine and $w^{*}$-$w^{*}$-continuous so that
$Af(y)=\phi(\bar{f}(\cdot + Ty))$ for all $y \in Y$.\\
$(ii)$ \ $\theta_{A}$ maps sheets into sheets.\\
$(iii)$  $\theta_{A}$ maps $Y^{**}$ into a sheet.

\noindent In particular, $\theta_A$ is continuous if and only if it is
continuous on $Y^{**}$
\end{theorem}
\begin{proof}%
Let us note first that $T:Y^{**} \longrightarrow X^{**}$ is affine
and $w^{*}$-$w^{*}$-continuous if and only if there exist $R:X^{*}
\longrightarrow Y^{*}$ linear and continuous and $x_{0}^{**} \in
X^{**}$ so that $T(y^{**})=R'(y^{**})+x_{0}^{**}$.

We begin by assuming that $(i)$ holds. If $A$ has such a
representation, let us see that then the Aron-Berner extension of
$Af$ is of the form
\begin{equation} \label{AronBerner}
\overline{Af}(y^{**})=\phi( \bar{f}( \cdot + Ty^{**})).
\end{equation}
Indeed, let $h(z)=\phi \big(f(\cdot + z) \big) = \phi \big(x \mapsto
f(x + z) \big)$ for $z \in X$. By \cite[Theorem 6.12]{ArCoGa91} its
Aron-Berner extension is given by $\bar{h}(x^{**})= \phi
\big(f(\cdot + x^{**}) \big) = \phi \big(x \mapsto f(x + x^{**})
\big)$.

We define 
$\tilde{h}(y^{**})= \phi \big(\bar{f}(\cdot +Ty^{**}) \big)$. Then
\[
\tilde{h}(y^{**})= (\bar{h} \circ T)(y^{**})=\bar{h} \big( R'(y^{**})+x_{0}^{**} \big)
= \big( \tau_{x_{0}^{**}}(\bar{h}) \circ R'\big) (y^{**}).
\]
Since $\bar{h}$ is the Aron-Berner extension of a function,
$\tau_{x_{0}^{**}}(\bar{h})$ is the Aron-Berner extension of some
other function (use, for example, \cite[Theorem 6.12]{ArCoGa91}). On
the other hand, by \cite[Lemma 9.1]{ArCoGa91} the composition of an
Aron-Berner extension with the transpose of a linear mapping is
again the Aron-Berner extension of some function. Hence
$\tilde{h}=\tau_{x_{0}^{**}}(\bar{h}) \circ R'$ is the Aron-Berner
extension of a function; but $\tilde{h}$ coincides with $Af$ on $X$,
therefore $\tilde{h}=\overline{Af}$ and \eqref{AronBerner} holds.

Now, to see that $\theta_{A}$ maps sheets into sheets it is enough
to find $S:Y^{**} \longrightarrow X^{**}$ such that
$\theta_{A}(\psi \circ \tau_{y^{**}}) = (\theta_{A} \psi) \circ
\tau_{Sy^{**}}$. We define $Sy^{**}=Ty^{**}+x_{0}^{**}$. First we
have
\begin{multline*}
\theta_{A} \big( \psi \circ \tau_{y^{**}} \big) (f) = \big( \psi \circ \tau_{y^{**}} \big)(A f) \\
= \psi [ y \mapsto \overline{A f} (y+y^{**})]
= \psi [ y \mapsto \phi [x \mapsto \bar{f} \big( x+ T(y+y^{**}) \big)]] \\
= \psi [ y \mapsto \phi [x \mapsto \bar{f} \big( x+ T y+ Sy^{**}
\big)]].
\end{multline*}
Let us call $g(x)= \bar{f}(x+Sy^{**})$. As above, we can check that
its Aron-Berner extension is $\bar{g}(x^{**})=
\bar{f}(x^{**}+Sy^{**})$. With this we obtain
\begin{eqnarray*}
\big( \theta_{A} \psi \circ \tau_{Sy^{**}} \big) (f)& =
&\theta_{A} \psi [ x \mapsto \bar{f}(x+Sy^{**})]
= \psi (Ag) \\
&=& \psi [ y \mapsto Ag(y) ] = \psi [ y \mapsto \phi [x \mapsto \bar{g}(x+Ty)]] \\
&=& \psi [ y \mapsto \phi [x \mapsto \bar{f} \big( x+ T y+ Sy^{**}
\big)]]
\end{eqnarray*}
and $(ii)$ holds. Clearly, $(ii)$  implies $(iii)$.

Let us suppose that $\theta_{A}$ maps $Y^{**}$ into a single sheet.
Hence, $\theta_{A}(\delta_{y^{**}}) = \theta (\delta_{0}) \circ
\tau_{Sy^{**}}= \phi \circ \tau_{Sy^{**}}$ for some $Sy^{**}$ in
$X^{**}$. This means that $\delta_{y^{**}}(Af) = \big( \phi \circ
\tau_{Sy^{**}} \big) (f)$ for all $f$ and from this
$\overline{Af}(y^{**})=\phi( \bar{f}( \cdot + Sy^{**}))$. Let us see
that $S$ is affine.

 Let $x^{*} \in
X^{*}$, then $Ax^{*}$ is a degree one polynomial and so is
$\overline{Ax^{*}}$. Also,
\begin{align*}
\overline{Ax^{*}}(y^{**}) & = \phi [ x \mapsto AB(x^{*}) (x+Sy^{**})] \\
& = \phi [ x \mapsto x^{*}(x) + Sy^{**}(x^{*})]
= \phi(x^{*}) + S(y^{**})(x^{*}).
\end{align*}
This shows that $S$ is $w^{*}$ affine; hence $S$ is affine.

 Let us
finish by proving that $S$ is $w^{*}$-$w^{*}$-continuous. Indeed, let
$(y^{**}_{\alpha})_\alpha$ be a net $w^{*}$-converging to $y^{**}$.
By  Lemma \ref{serafin} we have, for every $x^{*} \in X^{*}$,
$Ax^{*}=y^{*}_{x^{*}}+ \lambda_{x^{*}}$. Then
$\overline{Ax^{*}}(y_{\alpha}^{**})=y_{\alpha}^{**}(y^{*}_{x^{*}})+
\lambda_{x^{*}}$ and this converges to $y^{**}(y^{*}_{x^{*}})+
\lambda_{x^{*}}=\overline{Ax^{*}}(y^{**})$. Finally, $\lim_{\alpha}
S(y_{\alpha}^{**}) = \lim_{\alpha}
\overline{Ax^{*}}(y_{\alpha}^{**}) - \phi(x^{*})=
\overline{Ax^{*}}(y^{**}) - \phi(x^{*}) =S(y^{**})(x^{*})$ and this
completes the proof.
\end{proof}

The previous theorem characterizes the homomorphisms $A$ for which
$\theta_{A}$ maps $Y^{**}$ into a sheet. A particular case is when
$Y^{**}$ is mapped precisely to $X^{**}$. These are those for
which $\phi = \delta_{T_{1}(0)}$ for some $T_{1}$. Then
\[
\overline{Af}(y^{**}) = \delta_{T_{1}(0)} [ x \mapsto \bar{f}(x + Ty^{**})] = \bar{f} \big(T_{1}(0)+Ty^{**} \big)
= \big( f \circ T_{2} \big) (y^{**}).
\]

Following \cite{CaGaMa05}, we say that $A:HV(X)\rightarrow HV(Y)$ is an
$AB$-composition homomorphism if there exists $g:Y^{\ast \ast} \to
X^{\ast \ast }$ such that
$\overline{A(f)}(y^{\ast\ast})=\overline{f}(g(y^{\ast\ast}))$ for
all $f\in HV(X)$ and all $y^{\ast\ast}\in Y^{\ast\ast}$. By the
proof of the previous theorem, if $A$ is an $AB$-composition
homomorphism, then $g$ must be affine. We can state the following:

\begin{corollary} \label{abcomposicion}
Let $X$ and $Y$ be symmetrically regular Banach spaces and $A:HV(X)
\rightarrow HV(Y)$ an algebra homomorphism. Then
$\theta_{A}(Y^{**})\subset X^{**}$ if and only if $A$ is the
$AB$-composition homomorphism associated to an affine mapping.
\end{corollary}

We feel that some important differences between the weighted
algebras studied here and the algebra of holomorphic functions of
bounded type are worthy to be stressed. By Theorem~\ref{final} and
the comments following it, any $AB$-composition homomorphism
induces a continuous $\theta_A$. In \cite{CaGaMa05}, examples are
presented of composition homomorphisms inducing discontinuous
$\theta_A$. Also, there are examples of homomorphisms for which
the induced mapping $\theta_A$ is continuous on $Y^{**}$ but is
not continuous on the whole $\mathfrak{M}_b(Y)$ (i.e., splits some
sheet other than $Y^{**}$ into many sheets). Note that these
homomorphisms are associated to composition operators given by
polynomials of degree strictly greater than one, and would not
work for $HV(X)$.

\bigskip
A consequence of Corollary~\ref{abcomposicion} is that, unless the
spectrum of $HV(X)$ coincides with $X^{**}$, there are homomorphisms
on $HV(X)$ that are not $AB$-composition ones. Indeed, given any
$\psi \in \mathfrak{M}_b(X)$, we can proceed as in the proof of
Theorem~\ref{final} to obtain a homomorphism that maps $Y^{**}$ into
the sheet containing $\psi$. If $\psi$ does not belong to
$X^{**}$, the homomorphism thus obtained is not an $AB$-composition
one.

The one to one identification $\mathfrak{M}_{b}(X) \hookrightarrow
\esv$ leaves   $X^{**}$  invariant. If there exists a polynomial on
$X$ that is not weakly sequentially continuous, then
$\mathfrak{M}_{b}(X)$ properly contains $X^{**}$ and then so does
$\esv$. Therefore, if there are polynomials on $X$ that are not
weakly sequentially continuous, then there are homomorphisms on
$HV(X)$ other than $AB$-composition ones.


\begin{thebibliography}{10}

\bibitem{AlArDi84} {
R.~Alencar, R.~M. Aron, and S.~Dineen.}
\newblock A reflexive space of holomorphic functions in infinitely many
  variables.
\newblock {Proc. Amer. Math. Soc.}, 90(3) (1984), 407--411.

\bibitem{ArBe78}
{ R.~M. Aron and P.~D. Berner.}
\newblock {A Hahn-Banach extension theorem for analytic mappings.}
\newblock { Bull. Soc. Math. Fr.}, 106 (1978), 3--24.



\bibitem{ArCoGa91}
{ R.~M. Aron, B.~Cole, and T.~W. Gamelin.}
\newblock {Spectra of algebras of analytic functions on a Banach space.}
\newblock { J. Reine Angew. Math.}, 415 (1991), 51--93.

\bibitem{ArGaGaMa96}
{ R.~M. Aron, P.~Galindo, D.~Garc{\'{\i}}a, and M.~Maestre. }
\newblock Regularity and algebras of analytic functions in infinite dimensions.
\newblock { Trans. Amer. Math. Soc.}, 348(2) (1996), 543--559.

\bibitem{BiBoGa93}
{K.~D. Bierstedt, J.~Bonet, and A.~Galbis.}
\newblock {Weighted spaces of holomorphic functions on balanced domains}.
\newblock { Michigan Math.}, 40 (1993), 271--297.

\bibitem{BiBoTa98}
{K.~D. Bierstedt, J.~Bonet, and J.~Taskinen.}
\newblock {Associated weights and spaces of holomorphic functions}.
\newblock { Studia Math.}, 127(2) (1998), 137--168.

\bibitem{Bo38}
{R.~P. Boas, Jr.}
\newblock Representations for entire functions of exponential type.
\newblock { Ann. of Math. (2)}, 39(2) (1938), 269--286.

\bibitem{Bo54}
{R.~P. Boas, Jr.}
\newblock {Entire functions}.
\newblock Academic Press Inc., New York, 1954.

\bibitem{Bo03}
{J.~Bonet.}
\newblock {Weighted spaces of holomorphic functions and operators between
  them}.
\newblock In U.~d.~S. Secretariado~de Publicaciones, editor, { Proceedings
  of the seminar of Mathematical Analysis (Univ. M\'alaga, Univ. Sevilla)},
  pages 117--138, Sevilla, 2003.

\bibitem{BoDoLi99}
{J.~Bonet, P.~Doma\'nski, and M.~Lindstr\"{o}m.}
\newblock {Essential norm and weak compactness of composition operators on
  weighted spaces of analytic functions}.
\newblock { Canad. Math. Bull.}, 42(2) (1999), 139--148.

\bibitem{BoDoLiTa98}
{J.~Bonet, P.~Doma\'nski, M.~Lindstr\"{o}m, and J.~Taskinen.}
\newblock {Composition operators between wighted Banach spaces of analytic
  functions}.
\newblock { J. Austral. Math. Soc. (Series A)}, 64 (1998) 101--118.

\bibitem{BoFr02}
{J.~Bonet and M.~Friz.}
\newblock {Weakly compact composition operators on locally convex spaces}.
\newblock { Math. Nachr.} 245 (2002), 26--44.

\bibitem{BoFrJo03}
{J.~Bonet, M.~Friz, and E.~Jord\'a.}
\newblock {Composition operators between weighted inductive limits of sapces of
  holomorphic functions}.
\newblock { Publ. Math. Debrecen}, {\bf 67} no. 3-4  (2005),  333--348.

\bibitem{CaCaGa00}
{F.~Cabello~S{\'a}nchez, J.~M.~F. Castillo, and
R.~Garc{\'{\i}}a.}
\newblock Polynomials on dual-isomorphic spaces.
\newblock { Ark. Mat.}, 38(1)  (2000), 37--44.

\bibitem{CaGaMa05}
{D.~Carando, D.~Garc{\'{\i}}a, and M.~Maestre.}
\newblock Homomorphisms and composition operators on algebras of analytic
  functions of bounded type.
\newblock { Adv. Math.}, 197(2) (2005), 607--629.

\bibitem{CaLa04}
{D.~Carando and S.~Lassalle.}
\newblock {$E'$} and its relation with vector-valued functions on {$E$}.
\newblock { Ark. Mat.}, 42(2) (2004), 283--300.

\bibitem{DaGa89}
{A.~Davie and T.~Gamelin.}
\newblock {A theorem on polynomial-star approximation.}
\newblock { Proc. Amer. Math. Soc.}, 106(2) (1989), 351--356.

\bibitem{DiDi98}
{J.~C. D{\'{\i}}az and S.~Dineen.}
\newblock Polynomials on stable spaces.
\newblock { Ark. Mat.}, 36(1) (1998), 87--96.

\bibitem{LibroDi99}
{S.~Dineen.}
\newblock { {Complex analysis on infinite dimensional spaces}}.
\newblock Springer Verlag, London, 1999.

\bibitem{GalMaRu00}
{P.~Galindo, M.~Maestre, and P.~Rueda.}
\newblock Biduality in spaces of holomorphic functions.
\newblock { Math. Scand.}, 86(1) (2000), 5--16.

\bibitem{GaMaRu00}
{D.~Garc{\'{\i}}a, M.~Maestre, and P.~Rueda.}
\newblock Weighted spaces of holomorphic functions on {B}anach spaces.
\newblock { Studia Math.}, 138(1) (2000), 1--24.

\bibitem{GaMaSe04}
{D.~Garc\'{\i}a, M.~Maestre, and P.~Sevilla-Peris.}
\newblock {Composition operators between weighted spaces of holomorphic
  functions on Banach spaces}.
\newblock { Ann. Acad. Sci. Fenn. Math.}, 29 (2004), 81--98.

\bibitem{GaMaSe06}
{D.~Garc{\'{\i}}a, M.~Maestre, and P.~Sevilla-Peris.}
\newblock Weakly compact composition operators between weighted spaces.
\newblock { Note Mat.}, 25(1) (2005/06), 205--220.

\bibitem{GarrJa02}
{M.~I. Garrido and J.~A. Jaramillo.}
\newblock Variations on the {B}anach-{S}tone theorem.
\newblock { Extracta Math.}, 17(3) (2002), 351--383.
\newblock IV Course on Banach Spaces and Operators (Spanish) (Laredo, 2001).

\bibitem{Ka70}
{N.~J. Kalton.}
\newblock Schauder decompositions in locally convex spaces.
\newblock { Proc. Cambridge Philos. Soc.}, 68 (1970), 377--392.

\bibitem{LaZa00}
{S.~Lassalle and I.~Zalduendo.}
\newblock To what extent does the dual {B}anach space {$E'$} determine the
  polynomials over {$E$}?
\newblock { Ark. Mat.}, 38(2) (2000), 343--354.

\bibitem{Ou93}
{L.~Oubbi.}
\newblock Weighted algebras of continuous functions.
\newblock { Results Math.}, 24(3-4) (1993), 298--307.

\bibitem{Vi06}
{D.~M. Vieira.}
\newblock Theorems of {B}anach-{S}tone type for algebras of holomorphic
  functions on infinite dimensional spaces.
\newblock { Math. Proc. R. Ir. Acad.}, 106A(1):97--113, 2006.

\end{thebibliography}

\section*{Acknowledgements}
We would like to thank our friends: J.~Bonet for all the help solving the difficulties
with the associated weights, especially with the definite statement and proof of
Proposition \ref{caract alg} and showing to us Proposition \ref{eta asoc}, M.~Maestre
for Example \ref{manolo} and many discussions together with D.~Garc\'{\i}a
that helped to improve the final shape of the article. We would also like to thank K.~D.~Bierstedt for useful remarks and comments.

Most of the work in this article was done while the second cited author was
visiting the Department of Mathematics of the Universidad de Buenos Aires during
the summer/winter of 2006 supported by grants GV-AEST06/092 and UPV-PAID-00-06. He
wishes to thank all the people in and outside the Department that made that such a
delightful time.

\end{document}